\RequirePackage{fix-cm}
\documentclass[smallextended]{svjour3}       % onecolumn (second format)
\smartqed  % flush right qed marks, e.g. at end of proof
\usepackage{graphicx}
\usepackage{mathptmx}      % use Times fonts if available on your TeX system
\usepackage{latexsym}
\usepackage{amsmath}
\usepackage{amsfonts}
\usepackage{amssymb}
\usepackage{algorithm}
\usepackage[noend]{algpseudocode}
\usepackage{cite}
\RequirePackage[colorlinks,citecolor=blue,urlcolor=blue]{hyperref}
\newcommand {\bfy}  { {\bf y} }

\newcommand {\bfomega}  { {\boldsymbol \omega} }

\newcommand{\hf}{\frac12}

\newcommand{\grad}{\ensuremath {{\bf{ \nabla}}}}

\newcommand{\bfd}{{\bf d}}

\newcommand{\bfe}{{\bf e}}
\newcommand{\bft}{{\bf t}}
\newcommand{\bfs}{{\bf s}}
\newcommand{\bfg}{{\bf g}}
\newcommand{\bfv}{{\bf v}}
\newcommand{\bff}{{\bf f}}
\newcommand{\bfA}{{\bf A}}
\newcommand{\bfP}{{\bf P}}

\newcommand{\bfI}{{\bf I}}

\newcommand{\bfV}{{\bf V}}

\journalname{}
\begin{document}

\title{How To Catch A Lion In The Desert - On The Solution Of The Coverage Directed Generation (CDG) Problem%\thanks{Grants or other notes
}

\titlerunning{How To Catch A Lion In The Desert - On The Solution Of The CDG Problem}        % if too long for running head

\author{ Raviv Gal \and Eldad Haber \and Brian Irwin*\thanks{* Corresponding author.} \and Bilal Saleh \and Avi Ziv}


\institute{Raviv Gal \and Bilal Saleh \and Avi Ziv \href{https://orcid.org/0000-0002-6309-250X}{\includegraphics[scale=0.4]{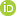}}
 	      \at
              IBM Research Laboratory in Haifa, Haifa, Israel \\
              % Tel.: +123-45-678910\\
              % Fax: +123-45-678910\\
              \email{\{RAVIVG, BILAL, AZIV\}@il.ibm.com}           %  \\
           \and
           Eldad Haber \and Brian Irwin \href{https://orcid.org/0000-0002-6086-4359}{\includegraphics[scale=0.4]{ORCIDiD_icon16x16.png}}
	      \at
              Department of Earth and Ocean Science, The University of British Columbia, Vancouver, BC, Canada \\
              \email{\{haber, birwin\}@eoas.ubc.ca}
}

\date{Received: date / Accepted: date}

\maketitle

\begin{abstract}
The testing and verification of a complex hardware or software system, such as modern integrated circuits (ICs) found in everything from smartphones to servers, can be a difficult process. One of the most difficult and time-consuming tasks a verification team faces is reaching coverage closure, or hitting all events in the coverage space. \emph{Coverage-directed-generation} (CDG), or the automatic generation of tests that can hit hard-to-hit coverage events, and thus provide coverage closure, holds the potential to save verification teams significant simulation resources and time. In this paper, we propose a new approach to the CDG problem by formulating the CDG problem as a noisy derivative free optimization (DFO) problem. However, this formulation is complicated by the fact that derivatives of the objective function are unavailable, and the objective function evaluations are corrupted by noise. We solve this noisy optimization problem by utilizing techniques from direct optimization coupled with a robust noise estimator, and by leveraging techniques from inverse problems to estimate the gradient of the noisy objective function. We demonstrate the efficiency and reliability of this new approach through numerical experiments with an abstract model of part of IBM's \emph{NorthStar} processor, a superscalar in-order processor designed for servers.
\end{abstract}

%==========================
\section{Introduction}
\label{sec:intro}
%==========================
Verification of a complex hardware or software system, such as modern \emph{integrated circuits} (ICs), can be a challenge. In principle, one would like to test every state or event that the system can reach, and observe that the system functions as intended. However, for complex systems, this is impossible, as the number of possible states is so large that it is impractical to test each state individually. To this end, it is common to define a large, but finite, random set of tests or \emph{test instances}, also referred to as \emph{test stimuli}, that are drawn from the distribution of all possible tests, and apply them to the \emph{design-under-test} (DUT) to be tested.

This paper targets verification environments that utilize \emph{biased random stimuli generators} to generate test stimuli. The stimuli generator uses, as its input, \emph{test templates} that bias the generation toward targeted areas and features of the verified design. A test template comprises a set of parameters, or directives, where each parameter is a set of weight-value pairs. We refer to the output stimuli of the random stimuli generator as a test instance.

Even with a smart choice of test stimuli, one may have great difficulty hitting a number of key events to be tested. These events are often referred to as \emph{hard-to-hit} events. This is because the mapping from test parameters to events is unknown, and can be highly nontrivial. Therefore, one of the most difficult and time-consuming tasks a verification team faces is reaching coverage closure, or, in other words, hitting all coverage events, including hard-to-hit events. Understanding why certain events are difficult to hit, and how they can be hit, requires both verification expertise and a deep understanding of the design under test. Moreover, generating test instances that hit such events is often an iterative trial and error process that consumes significant simulation resources and verification team time. Therefore, it is desirable to have an  automatic solution for improving the probability of hitting hard-to-hit events.

\emph{Coverage-directed-generation} (CDG), or the automatic generation of test instances, is a concept that has long been on the wish list of verification teams, and the target of a vast amount of research. Many techniques have been proposed to tackle the CDG problem, ranging from formal methods, via AI algorithms, to data analytics and machine learning techniques (see \cite{pipe-prev4}, \cite{gilly_cdg}, \cite{CDGFunctionalVerificationBayesianNetworks} and references within). These techniques did not mature to be widely used in industry for various reasons, including the scalability of the solution, difficulty in applying it, and the quality of the proposed solution. As a result, reaching coverage closure remains almost entirely a manual process.

The goal of this work is to propose a new approach for the solution of the CDG problem and increasing the probability of hitting hard-to-hit events. Finding how to hit a low probability event in a large space is sometimes humourously referred to as finding ``how to catch a lion in the  desert", originated by the seminal paper of P\'{e}tard \cite{HowToCatchALionPetard}. We propose a method that solves the problem by minimizing a cost function that increases the probability of hitting the hard-to-hit event(s). We show that such an approach can lead to an efficient solution of the problem, especially if it is coupled with a robust and efficient optimization algorithm.

The rest of the paper is organized as follows. In Section~\ref{sec:mathematical-background-sec}, we give a mathematical background to the proposed approach. In Section~\ref{sec:solution-technique-sec}, we discuss solution techniques for the problem. These techniques are based on direct optimization methods coupled with a robust noise estimator. In Section~\ref{sec:experimental-setup-NorthStar}, we describe the main experimental environment used to test the proposed approach. In Section~\ref{sec:numerical-experiments-sec}, we perform a number of experiments that demonstrate the efficiency of our approach, and we summarize the paper in Section~\ref{sec:conclusions-sec}.

%==========================
\section{Mathematical Background}
\label{sec:mathematical-background-sec}
%==========================
Let us mathematically formalize the testing process. Let $\theta(\bft)$ denote a random variable, referred to as a {\em test instance} of {\em test template} $\bft$, and representing a test to be run by the DUT. Using {\em directives}, the test template $\bft$ can be represented as a vector $\bft = [\bfd_1,\ldots,\bfd_n]$ composed of $n$ {\em directive weight vectors} $\bfd_j, j=1,\ldots,n$. Each directive $d_j$ is parametrized by a weight vector $\bfd_j$, where the weight vector is normalized to present a probability distribution. The space of all possible test templates is denoted by ${\cal T}$, and is also known as the {\em test templates skeleton}. It is important to note that, while each test instance $\theta(\bft)$ is random, the directives and the test templates are \textbf{not}. The directives and test templates are deterministic parameters that define the random space and control the distribution of the test instances. 

In the testing and verification process, the main goal is to hit every event in the {\em coverage space} ${\cal C} = \{c_1, ..., c_m\}$, or space of all events. Given a test instance $\theta(\bft)$, chosen from a probability space defined by the vector $\bft \in {\cal T}$, one runs a simulation to obtain a random vector, $\bfs(\theta) = [s_1,\ldots, s_m], \ s_j \in \{0,1\} \text{ } \forall j$, that is defined as a {\em hit coverage} vector. The entries of the hit coverage vector $\bfs$ are binary. If a particular event in the coverage space was hit by the specific test instance $\theta$, the entry of the corresponding index in $\bfs$ is $1$, and it is $0$ otherwise.

Clearly, since the test instances are generated randomly in a manner dependent on the parameters of the test template $\bft$, the vector $\bfs$ is also random and depends on $\bft$. To this end, let
\begin{eqnarray}
\label{eqn:hit-coverage}
{\bfe}(\bft) = {\mathbb E} \, \big [ \bfs ({\theta(\bft)}) \big ]
\end{eqnarray}
be the expected value of the hit coverage vector $\bfs$, and let
\begin{eqnarray}
\label{eqn:empirical-hit-coverage}
\bfe_N(\bft) = {\frac 1N} \sum_j \bfs_j
\end{eqnarray}
be the empirical expectation of the hit coverage vector estimated using $N$ test instances generated from test template $\bft$. Note that while $s_j \in \{0,1\} \text{ } \forall j$, the vector ${\bfe} = [e_1,\ldots,e_m]$ and its empirical values are {\em real}. The $j$-th value in $\bfe(\bft)$, $e_j$, represents the probability of hitting the event $c_j$ using a test instance $\theta$ generated according to the distribution defined by $\bft$. To compute the empirical expectation $\bfe_N(\bft)$ of a hit coverage vector, given a test template $\bft$, we can run $N$ simulations, obtain $\bfs_j, j=1,\ldots,N$ hit coverage vectors, and average them. Clearly, such a process is computationally expensive, especially if we are to estimate $\bfe(\bft)$ accurately. To demonstrate the above definitions, let us consider the following simple, but concrete example.

\paragraph{\bf Example 1: Testing The Multiplication of Two Numbers}
Assume that we build a calculator that can compute the product of two numbers in the interval $[0,1]$. In a test instance, we need to randomly pick two numbers within the interval and compute their product.
In this case, we have two directive weight vectors, $\bfd_1$ and $\bfd_2$, that define how we choose each of the two numbers. For simplicity, in this case we assume that $\bfd_1=\bfd_2 = \bfd$, and therefore the test template $\bft$ is just the single directive weight vector $\bft = \bfd$. Next, we choose  the parametrization of the test template, which defines the numbers in the interval $[0,1]$. For simplicity, we assume that $\bft = [t_1,\ldots,t_k]$ are the probabilities of choosing a number in the interval $[0,1/k],\ldots,(1-1/k,1]$. 

Recall the space ${\cal C}$ is the space of all events, or coverage space. Let us define $m=k$ different events that correspond to the $k$ cases that the output of the multiplication falls into one of the intervals $[0,1/k],\ldots,(1-1/k,1]$. Now, consider choosing the probability density parameterized by $\bft$. One tempting choice is to simply use the uniform distribution, setting $t_i = 1/k, i=1,\ldots, k$. Clearly, this choice leads to less than optimal sampling of the coverage space. For this case, it is easy see that 
\begin{equation*}
e_1 \gg e_k
\end{equation*}
If we further refine the intervals in the ${\cal T}$ and ${\cal C}$ spaces by letting $k \rightarrow \infty$, then the likelihood of hitting an event that is at the right edge (close to $1$) will approach $0$, and therefore using a uniform distribution may not lead to a complete sampling of the coverage space, and we may end up with some unhit events. Understanding this problem allows one to choose a sampling routine that gives a higher probability to numbers that are close to $1$, and improve the probability of sampling the whole coverage space. \\

The above multiplication example can be clearly analyzed to obtain an optimal sampling scheme. However, in practice, this is very frequently not the case. The system under test may be highly nonlinear. In this case, one typically performs some probing of the space by randomly testing a number of sampling schemes, and then tries to improve the coverage and sample as intelligently as possible. However, as previously discussed, hitting a hard-to-hit event may be difficult and require manual and labor intensive processes. Our goal is to improve over such processes by {\em automatically} increasing the probability of hitting hard-to-hit events.

Obtaining $\bfe_N(\bft)$ from $\bft$ is an unknown function that is dictated by the simulator, and can be written  as
\begin{eqnarray}
\label{eqn:noisy-empirical-hit-coverage}
\bfe_N(\bft) = \bfe(\bft)  + \bfomega(\bft)
\end{eqnarray}
Here, $\bfomega(\bft)$ is a noise vector that depends on the parameters, $\bft$. This noise vector $\bfomega$ gets a different value every time we compute $\bfe_N$, giving us a noisy realization of the expected value of the hit coverage vector. 

Let us define the target event(s), $\bfe^{\rm tar}(\bft) = \bfP^{\top}  \bfe(\bft)$. Depending on the specific problem, $\bfe^{\rm tar}$ can be a vector or a number. For example, if we only want to hit the $j$-th event, we can define $\bfP^{\top} $ as the $j$-th row of the identity matrix. In some cases, we aim to increase the probability of hitting a group of hard-to-hit events, and in this case $\bfP^{\top}$ corresponds to a few rows of the identity matrix. Maximizing the probability of hitting the events in $\bfe^{\rm tar}$ can now be formulated as the simple optimization problem
\begin{eqnarray}
\label{eqn:simple-coverage-optimization}
\max_{\bft}  \bigg \{ \phi(\bft) = {\boldsymbol 1}^{\top}\bfe^{\rm tar}(\bft) =   {\mathbb E} \big [  {\boldsymbol 1}^{\top}\bfP^{\top} \bfs(\theta(\bft)) \big ] \bigg \}
\end{eqnarray}
where ${\boldsymbol 1}$ is a vector of all ones.

There are a number of problems when attempting to solve the maximization problem defined by equation~\ref{eqn:simple-coverage-optimization}. First, we do not have access to the objective function directly. The objective function can only be evaluated up to some unknown noise. Second, this noise is not necessarily stationary. That is, every time the objective function is called, a different noise vector is generated, and, on top of that, the noise level can be different for different values of $\bft$. Third, for a fixed $\bft$, the noise corrupting the measurement of $\bfe(\bft)$ is likely different for each entry. In other words, the noise level is likely different for each event $e_j$. Fourth, a critical difference between this problem and the common problem of minimization under the expectation is that for the canonical stochastic programming problem, the random variable is drawn from a {\em fixed} distribution. Here, the distribution is parameterized by $\bft$, and therefore, as we change the values of the parameters we optimize, we obtain a different distribution with a possibly different noise signature. To illustrate the above, we continue with our discussion of Example~1.\\

\begin{figure}[htbp]
\centering
\includegraphics[width=4.5in]{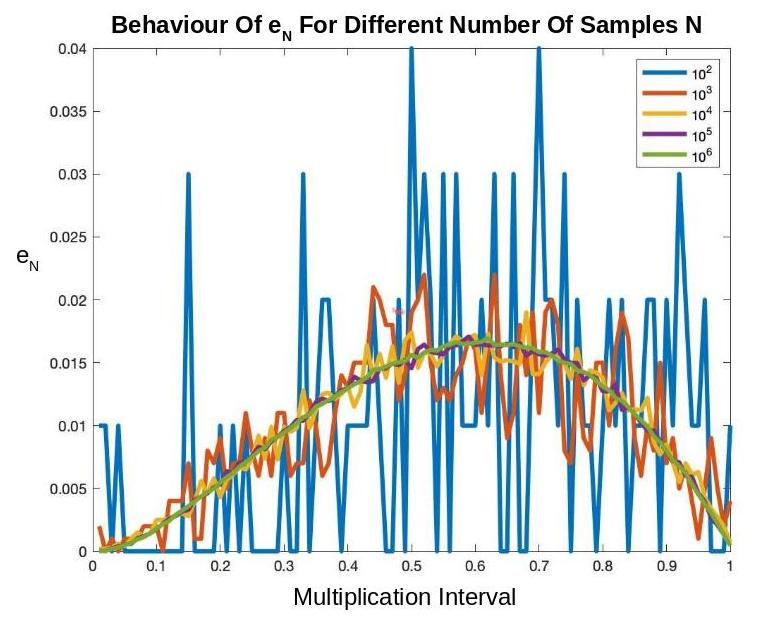}
\caption{Evaluation of $\bfe_N$ for the values $N = 10^i, i=\{2,3,4,5,6\}$ on the simple two number multiplication model problem. We choose $k=100$ segments, and choose the elements of $\bft$ to grow quadratically in the interval $[0,1]$, and normalize such that they sum to $1$. Note how noisy the function can be when the number of realizations is small, and how the probability converges as the number of samples grows. Also note the low probability of hitting events at either end (close to 0 or 1) of the multiplication interval.}
\label{fig:multiplication-simulator-sampling}
\end{figure}

\paragraph{\bf Example 1: Testing The Multiplication of Two Numbers - Continued} 
We choose $k=100$ segments, and choose the entries of $\bft$ to grow quadratically in the interval $[0,1]$, and normalize such that they sum to $1$.
This implies that we have a higher probability of choosing larger numbers compared with smaller numbers.
Given this test template, we compute the empirical hit coverage vector, $\bfe_N(\bft)$, for $N = 10^i, i=\{2,3,4,5,6\}$.
The results are plotted in Figure~\ref{fig:multiplication-simulator-sampling}.
The results demonstrate how noisy the function can be when the number of realizations is small, and how the probability converges as the number of samples grows. Also note how low the probability of choosing numbers close to 1 is, even when $\bft$ is chosen to grow quadratically. To find a $\bft$ that further improves the probability of hitting the rightmost element in the empirical hit coverage vector $\bfe_N$, by setting ${\boldsymbol 1}^{\top}\bfP^{\top} = [0,\ldots,0,1]$, one can compute an objective function that maximizes the probability of hitting the rightmost element.

%==========================
\section{Solution Techniques}
\label{sec:solution-technique-sec}
%==========================
The main problem under consideration here, the CDG problem, can be formulated as a derivative free optimization (DFO) problem where the objective function under consideration is noisy. The topic has been considered by many authors, such as C.T. Kelley, J. Nocedal, and K. Scheinberg, using different techniques, ranging from stochastic methods \cite{IntroductionToDFO}, direct search methods \cite{kelley3}, and gradient based methods \cite{DFONoisyFunctionsQuasiNewton}. In this paper, we experiment with three optimization techniques: an implicit filtering based technique, a steepest descent based technique, and a Broyden-Fletcher-Goldfarb-Shanno (BFGS) based technique. While we have directly used an implicit filtering technique, we have modified existing steepest descent and BFGS techniques from non-noisy unconstrained optimization in order to better deal with the noise in the problem. Below, we describe the algorithmic framework we used to solve the problem.

\subsection{Optimization Problem Setup}
\label{sec:optimization-problem-setup}
Given an objective function $f(\bft)$, we rewrite it as
\begin{eqnarray}
\label{eqn:objective-decomposition}
f(\bft) = \phi(\bft) + \omega(\bft)
\end{eqnarray}
We assume that $\phi(\bft)$ is a smooth function, and that $\omega(\bft)$ is noise. The noise $\omega$ is assumed to be uncorrelated with zero mean, and some unknown standard deviation $\sigma$. We assume that the standard deviation $\sigma(\bft)$ changes slowly with respect to $\bft$.

For the CDG problem, we are unable to obtain the derivatives of $\phi$ with respect to $\bft$, and therefore we turn to DFO methods. While there are many DFO methods, we turn our attention to local methods that are based on the numerical estimation of the gradient. Such methods have been studied extensively in the last 30 years \cite{DFOAlgReviewSoftwareComparison}, yielding successful software  packages such as MCS, TOMLAB/LGO, and NEWUOA \cite{GlobalOptimizationMCS, GlobalOptInAction, NEWUOA} (also see references within). 

However, when experimenting with the problem, we found that standard approaches based on gradient estimation methods fail or work poorly when the noise level is high. To explain this problematic observation, we first review the standard approach to such problems. A typical algorithm for such problems is composed of the following steps.
\begin{enumerate}
\item Evaluate the function $f$, its gradient $\nabla f$, and its approximate Hessian $\bf B$.
\item Compute a descent direction $\bf z$.
\item Update the solution using some relaxed line search or trust region method.
\end{enumerate}
Function and gradient evaluations are typically done using finite
differences. Let us review the process at some depth. Assume that we would like to compute the directional derivative
of $f(\bft)$ in the direction $\bfv$. It is common to use a central finite difference approach computing
\begin{eqnarray}
\label{eqn:directional-derivative-CFD-approx}
\bfv^{\top}\grad f(\bft) \approx  {\frac {f(\bft + h \bfv) - f(\bft - h \bfv)}{2h}}
\end{eqnarray}
which gives
\begin{eqnarray}
\label{eqn:directional-derivative-objective-decomposition-approx}
\bfv^{\top}\grad f(\bft) \approx \bfv^{\top} \grad \phi(\bft) + {\frac {\bar\omega}{2h}} + h^2N(\phi(\bft),\bfv)
\end{eqnarray}
where $\bar \omega$ is a random variable generated by combining the zero mean errors in the function evaluations. It is evident that the approximation for  $\bfv^{\top} \grad f(\bft)$  is polluted with two types of errors. The first type of error, corresponding to the second term in equation~\ref{eqn:directional-derivative-objective-decomposition-approx}, is the error due to the noisy estimation of the function, and the second type of error, corresponding to the third term in equation~\ref{eqn:directional-derivative-objective-decomposition-approx}, is an error due to the nonlinearity of $\phi(\bft)$. Unfortunately, these error terms have contradicting behaviours. While the second term in equation~\ref{eqn:directional-derivative-objective-decomposition-approx} requires as large an $h$ as possible to reduce the error, the third requires a small $h$ to obtain the same goal. In some cases, when the noise is small, and it is possible to obtain an estimate of the magnitude of the nonlinear residual, one can balance these terms, choosing
$$ h = \left({\frac {\sigma}{2 \bar N}} \right)^{\frac 13} $$ 
where $\bar N$ is an estimate of the nonlinear residual. This approximation can be used in order to
obtain a reasonable estimate of the gradient. However, even with this choice of approximation, the estimate of the gradient may not be sufficiently accurate.

 Furthermore, estimating the noise and the nonlinear errors can be computationally difficult, and require additional function evaluations using different sized stencils. Such work was proposed in \cite{MoreWildNoise}. However, even with an optimal stencil size, the noise can still be significant (see, for example, Figure~\ref{fig:multiplication-simulator-sampling}). Indeed, even for the optimal $h$ (assuming that both $\bar N$ and $\sigma$ are known), the error corrupting $\bfv^{\top}\grad f(\bft)$ scales as $\sigma^{\frac 23}$, which only marginally improves the problem presented by the noise for large values of $\sigma$.

In this work, we introduce a different approach to the optimization problem. Rather than estimating the noise by further function evaluations, we view the
problem as a statistical inverse problem, where the solution has to be evaluated from noisy data. In the next subsection, we show how to use standard techniques from inverse problems to estimate the behaviour of the objective function $f$, and its gradient $\nabla f$.
 
\subsection{Function And Gradient Approximation As Statistical Parameter Estimation}
\label{sec:approx-by-statistical-parameter-estimation}
Let us provide a different interpretation of the process of evaluating the gradient of a noisy function. Let us consider a general linear model of the form
\begin{eqnarray}
\label{eqn:f-general-linear-model}
f(\bft + h \bfv_i) =  \bar\phi(\bft) + h\bfv_i^{\top} \bfg + \omega_i 
\end{eqnarray}
with $\|\bfv_i\|=1$ and $\omega_i$ being noise. Here, $\bfg$ is an unknown vector that is to be computed from the values of the objective function in points around $\bft$. Note that this linear approximation is {\bf not} necessarily the Taylor expansion. It can be any linear model that  approximates the function for a given step size $h$ and direction $\bfv_i$. Clearly, for smooth functions, as $h\rightarrow 0$, the approximation converges to the Taylor expansion.

Now, assume that we have $n$ directions, $\bfv_1, \ldots, \bfv_n$. Using these directions, we obtain the following set of equations
\begin{eqnarray}
\label{eqn:mt}
\begin{pmatrix} f_1 \\  \vdots \\ f_n \end{pmatrix}
= \begin{pmatrix} 1 &  &h\bfv_1^{\top} & \\ 1 &  & \hdots & \\ 1 &  & h\bfv_n^{\top} &  \end{pmatrix} \begin{pmatrix} \bar \phi  \\ \bfg \end{pmatrix} + 
\begin{pmatrix} \omega_1  \\ \vdots \\ \omega_n \end{pmatrix}
\end{eqnarray}
which we rewrite as the simple linear system
\begin{eqnarray}
\label{eqn:mt1}
\bff =  \bfV \widehat \bfg + \boldsymbol \omega.
\end{eqnarray} 
where $\widehat \bfg = [\bar \phi, \bfg^{\top}]^{\top}$.
Estimating $\widehat \bfg$ from the noisy data $\bff$ is a corner stone of statistical inverse problems \cite{TenorioIntroDAandUQforIP}. It is therefore straight forward to use inverse problems techniques for the estimation of the average function value, $\bar \phi$, and gradient $\bfg$.

We further assume that we have some prior estimate of $\widehat \bfg$, $\widehat \bfg_0$. If no such estimate is available, then we can choose $\widehat \bfg_0 = \bold{0}$. Such an estimate can be obtained if we know something about the function $f$, or if we computed $\widehat \bfg$ at a nearby point. For example, if $\widehat \bfg$ was computed during a previous iteration, we can use this value from the previous iteration  as $\widehat \bfg_0$. A new estimate of $\widehat \bfg$ can be obtained by solving the following ridge regression type minimization problem
\begin{eqnarray}
\label{eqn:regression-objective}
\min_{\widehat \bfg}\ \bigg \{ \hf \|\bfV \widehat \bfg - \bff\|^2 + {\frac {\alpha}2} \|\widehat \bfg - \widehat \bfg_0\|^2 \bigg \}
\end{eqnarray}
Given the regularization parameter, $\alpha$, the problem has the closed form solution
\begin{eqnarray}
\label{eqn:regression-closed-form-solution}
\widehat \bfg_{\alpha} = (\bfV^{\top}\bfV + \alpha \bfI)^{-1}(\bfV^{\top} \bff + \alpha  \widehat \bfg_0)
\end{eqnarray}
The regularization parameter $\alpha$ is chosen based on the noise level. When the noise level is unknown, as in our problem, the Generalized Cross Validation (GCV) method can be used to choose $\alpha$, and obtain an unbiased estimate of the noise level \cite{ghw}. This is done by minimizing the GCV function for this problem
\begin{eqnarray}
\label{eqn:regression-GCV-function}
{\rm GCV}(\alpha) = {\frac {\| (\bfI - \bfA(\alpha))(\bff-\bfV \widehat \bfg_0)\|^2}{{\rm trace}\left( \bfI - \bfA(\alpha)\right)^2}}
\end{eqnarray}
where 
$$\bfA(\alpha) = \bfV(\bfV^{\top}\bfV + \alpha \bfI)^{-1}\bfV^{\top}. $$
Minimizing the GCV function in equation~\ref{eqn:regression-GCV-function} in 1D can be done using  a bisection method \cite{BurdenFairesNumericalAnalysis}. As we will later see in the numerical experiments section, the values obtained using the above approach can provide a significant advantage compared to the simple finite difference approximations employed in classical noisy optimization approaches. 

\subsection{Solution Algorithm}
\label{sec:solution-algorithms}

Below, we present our solution algorithm in pseudocode.  Algorithm~\ref{alg:gradient-based-steepest-descent} outlines the gradient based steepest descent technique and the BFGS approximation.

As a comparison to our approach, we use implicit filtering as presented in \cite{kelley3}. The implicit filtering algorithm does not require any gradient approximations, and only relies on the noisy values of $f$ itself. Our algorithm adapts standard descent methods for non-noisy problems by using the GCV estimated $\bar \phi$ and $\bfg$ in place of the values of $f$ and $\nabla f$ at each point, and performs a simple line search procedure.

A few comments are in order
\begin{itemize}
\item Minimizing the GCV function is not, in general, a computationally cheap process.
However, for the problem at hand, function evaluation is very expensive and the number of variables does not exceed the few thousands. In this case, investing some work to obtain
the best direction possible is justified.
\item For problems where function evaluation is cheap, one may not find our approach attractive.
\item Efficient ways to minimize the GCV function that use stochastic trace estimators can make the process of solving the problem relatively fast. Here we have used the technique proposed in \cite{GCVfast} to obtain the solution of the problem using Krylov space decomposition.
\end{itemize}

\begin{algorithm}[htbp]
\caption{Gradient Based Steepest Descent Technique}
\label{alg:gradient-based-steepest-descent}
\begin{algorithmic}
\Procedure{Descent Technique}{}
\State \% Initialize minimization algorithm
\State $iter \gets 0$ \% Iteration count
\State $\mu_{ls} \gets 10$ \% Initialize line search parameter
\State $ls_{break} \gets$ False \% Line search break flag
\State $\bft \gets \bft_{init}$
\State $\bft_{opt} \gets \bft_{init}$
\State $\bar \phi_{opt} \gets \infty$
\State Evaluate $f(\bft + h \bfv)$ in $n$ random directions $\bfv$
\State Estimate $\bar \phi$ and $\bfg$ by solving~\ref{eqn:regression-closed-form-solution} using GCV
\State Approximate the inverse Hessian, $\bf B$  (in the case of BFGS) or set $\bf B = \bf I$
  \While {True}
    \State $iter \gets iter + 1$ 
    \State $\bft_{old} \gets \bft$
    \State $\bar \phi_{old} \gets \bar \phi$
    \State $\bfg_{old} \gets \bfg$
    \State $lsIter \gets 1$ \% Line search iteration count
    \While {True}
      \State $\bft \gets \bft_{old} - \mu_{ls} {\bf B} \bfg_{old}$
      \State Evaluate $f(\bft + h \bfv)$ in $n$ directions $\bfv$
      \State Estimate $\bar \phi$, $\bfg$, and average noise level $||\omega||$
      \If {$\bar \phi < \bar \phi_{old} + 2||\omega||$}
        \If {$\bar \phi < \bar \phi_{opt}$}
          \State $\bft_{opt} \gets \bft$ 
          \State $\bar \phi_{opt} \gets \bar \phi$
          \State \textbf{break}
        \EndIf
      \EndIf
      \State $\mu_{ls} \gets \mu_{ls}/2$ \% Shrink line search parameter
      \State $lsIter \gets lsIter + 1$ 
      \If {$lsIter >$ Max. \# line search iterations}
        \State $ls_{break} \gets$ True \% Line search break
        \State \textbf{break}
      \EndIf
    \EndWhile
    \If {$lsIter = 1$}
      \State $\mu_{ls} \gets 2  \mu_{ls}$ \% Expand line search parameter
    \EndIf
    \State \% Check algorithm termination conditions
    \If {$ls_{break} =$ True} 
      \State \textbf{break}
    \EndIf
    \If {$iter >$ Max. \# iterations}
      \State \textbf{break} 
    \EndIf
  \EndWhile
\EndProcedure
\end{algorithmic}
\end{algorithm}

%%%

\bigskip

Before presenting the results from numerical experiments using our approach, we first describe the system used for the numerical experiments: an abstract model of part of IBM's NorthStar processor.

%==========================
\section{The NorthStar Pipeline}
\label{sec:experimental-setup-NorthStar}
%==========================
As a lightweight experimental environment, we employ a high-level software model of the two arithmetic pipes of the NorthStar superscalar in-order processor and the dispatch unit, also used in \cite{CDGFunctionalVerificationBayesianNetworks}. The NorthStar processor, also known as the RS64-II or PowerPC A50, was released by IBM in the late 1990s, featuring a RISC instruction set architecture \cite{Salvatore_4thgeneration}. The high-level software model consists of two main components. First, a biased random stimuli generator that generates programs, and second, a software simulator of the NorthStar processor's dispatch unit and two arithmetic pipes that executes the randomly generated programs.

\begin{figure}
  \centering
  \includegraphics[width=4.5in]{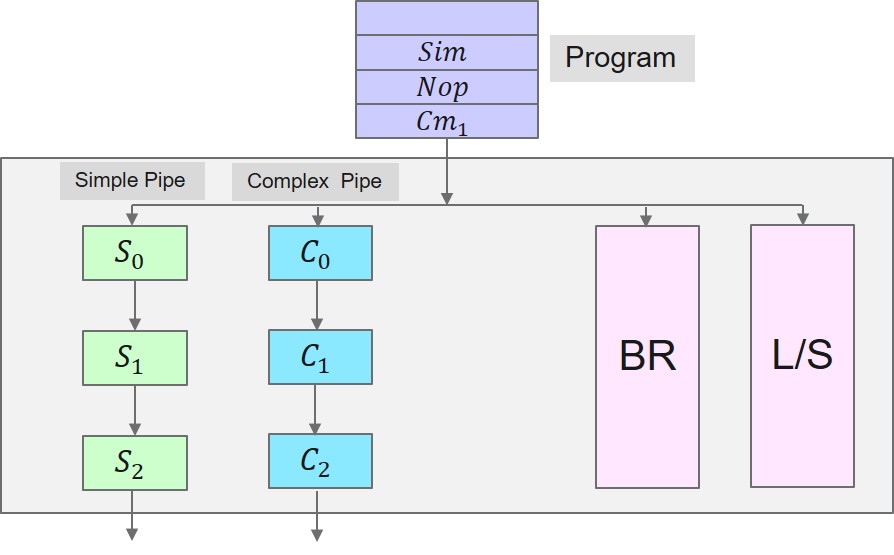}
  \caption{Schematic of the simulated NorthStar pipeline. There are two pipes of 3 stages, one simple pipe $S$ and one complex pipe $C$. In addition, L/S represents the processor's load store unit, and BR the branch prediction unit.}
  \label{fig:NorthStar-schematic}
\end{figure}

The NorthStar has two pipes, one simple and one complex (see Figure \ref{fig:NorthStar-schematic}). Each of the pipes comprises three stages: data fetch, execution, and write-back. One of the pipes, the simple pipe, handles only simple instructions, such as $add$. The other pipe, the complex pipe, handles complex instructions, such as $mul$. The complex pipe can also handle simple instructions when the simple pipe is busy. The model supports five types of instructions: simple instructions $Sim$, three types of complex instructions $Cm_1, Cm_2, Cm_3$ that differ in the time they spend in the execution stage (1, 2, and 3 cycles respectively), and $Nop$, which represents all instructions that are not executed in the arithmetic pipes. The actual execution time can be longer due to data dependencies between instructions. To maintain simplicity, we assume the processor has only eight registers and instructions use one source and one target register. In addition, the processor has a condition register $CR$, which some instructions read from and write to. In each cycle, up to two instructions are fetched, according to the instruction's type and the state of the pipes. 

Test templates $\bft$ for the NorthStar software model are defined by four directive weight vectors $\bft = [\bfd_1,\bfd_2,\bfd_3,\bfd_4]$, and control the distribution that the biased random stimuli generator component generates random programs from. Each directive weight vector defines a probability distribution. The first directive weight vector $\bfd_1 = IW = [W_{Nop}, W_{Sim}, W_{Cm_1}, W_{Cm_2}, W_{Cm_3}]$ contains instruction set selection weights, and controls the mnemonic of the generated instructions. The second and third directives affect the behaviour of the source and target registers. The second directive weight vector $\bfd_2 = SW = [W_{S_0}, ..., W_{S_7}]$ contains source register weights, and the third directive weight vector $\bfd_3 = TW = [W_{T_0}, ..., W_{T_7}]$ contains target register weights. The fourth directive weight vector $\bfd_4 = CW = [W_{C_0}, W_{C_1}]$ controls the conditional register. Thus, one can express a test template as the 23 entry vector $\bft = [IW, SW, TW, CW]$. 

% explain the $(Sim, nop, *, *, *)$
The coverage space $\cal C$ is a cross-product \cite{cov_book} of the instructions in stage 0 of the complex and simple pipes (5 and 2 possible values respectively), two indicators for whether stage 1 of each pipe is occupied, and an indicator for whether the instruction in $S_1$ is using the conditional register. An event is defined by assigning values to each coordinate. For example, the event $(C_2, Sim, 0, 0, 0)$ means that $C_2$ and $Sim$ are hosted at stage 0 of the complex and simple pipes, stage 1 of both pipes is not occupied, and the conditional register is not used. Clearly, the size of the coverage space size is $|C_0Inst \times S_0Inst \times C_1Used \times S_1Used \times S_1CR| = |5 \times 2 \times 2 \times 2 \times 2| = 80$. However, out of this space, only 54 events are legal. For example, the 8 events spanned by the subspace $(Sim, Nop, *, *, *)$, where * indicates a wildcard that can be any value, are illegal because if $S_0$ is free, then the simple instruction should have been fetched into the simple pipe. During simulation, coverage is tracked for a time interval of 100 cycles, starting at cycle 10. An event is considered hit by the test instance if it was hit at least once during this time interval.

%==========================
\section{Numerical Experiments}
\label{sec:numerical-experiments-sec}
%==========================
In this section, we compare the performance of the implicit filtering, steepest descent, and BFGS techniques numerically using the NorthStar pipeline simulator described above in section~\ref{sec:experimental-setup-NorthStar}.

\subsection{Initial Exploration}
\label{sec:initial-exploration}
As an initial exploration of $\bfe(\bft)$ for the NorthStar, we first ran 5000 random test templates drawn from ${\cal T}$ according to the Dirichlet distribution $Dir(1)$. Using these 5000 test templates, we hit all events in the coverage space $\cal C$ at least once. We also found the hardest event to hit to be event $c_{hard} = (C_2, Nop, 0, 1, 0)$. The single best test template hit event $c_{hard}$ with probability $p(c_{hard})=0.15$. Based on applied domain knowledge, it was deduced that the test template defined by $IW=(0.5, 0.2, 0, 0.3, 0)$, $SW=TW=(1, 0, 0, 0, 0, 0, 0, 0)$, and $CR=(1,0)$ would yield the best chance of hitting event $c_{hard}$\footnote{There are many other templates with different values of $SW$ and $TW$ that achieve the same probability.}. This test template 
\begin{equation*}
\bft = [0.5,0.2,0,0.3,0,1,0,0,0,0,0,0,0,1,0,0,0,0,0,0,0,1,0] 
\end{equation*}
will give high weights to $Nop$, $Sim$, and $Cm_2$, will create dependencies between the source and target registers, and will not use the condition register $CR$. Experimentally, by averaging over $100000$ runs of the simulator, this template was observed to yield a hit probability of $p(c_{hard})=0.4$. Below, we continue to use values averaged over $100000$ runs of the simulator to define a high quality estimate, or ``true" value of $p(c_{hard})$. However, as we show below, both the implicit filtering and steepest descent based techniques are able to {\em automatically} discover test templates that achieve $p(c_{hard})=0.4$ or close to $0.4$ with a modest budget of total runs of the NorthStar simulator. 

\subsection{Event $c_{hard}$ Objective Function}
\label{sec:c_hard-objective-function}
As was the case when analyzing how to maximize the probability of hitting the rightmost element in the empirical hit coverage vector for the two number multiplication simulator in Section~\ref{sec:mathematical-background-sec}, we can again choose $\bfP^{\top} $ to be a single row of the identity matrix. Specifically, $\bfP^{\top} $ is now the row of the identity matrix corresponding to $c_{hard}$. To avoid explicity enforcing the constraint that the directive weight vectors $IW, SW, TW, CW$ define properly normalized probability distributions, and thus solving a constrained optimization problem, we instead pass the values obtained from the optimization algorithms through the standard softmax function to ensure valid probability distributions before passing them to the program generator component of the NorthStar software model. 

\begin{figure}
  \centering
  \begin{tabular}{cc}
  \includegraphics[width=5.1cm]{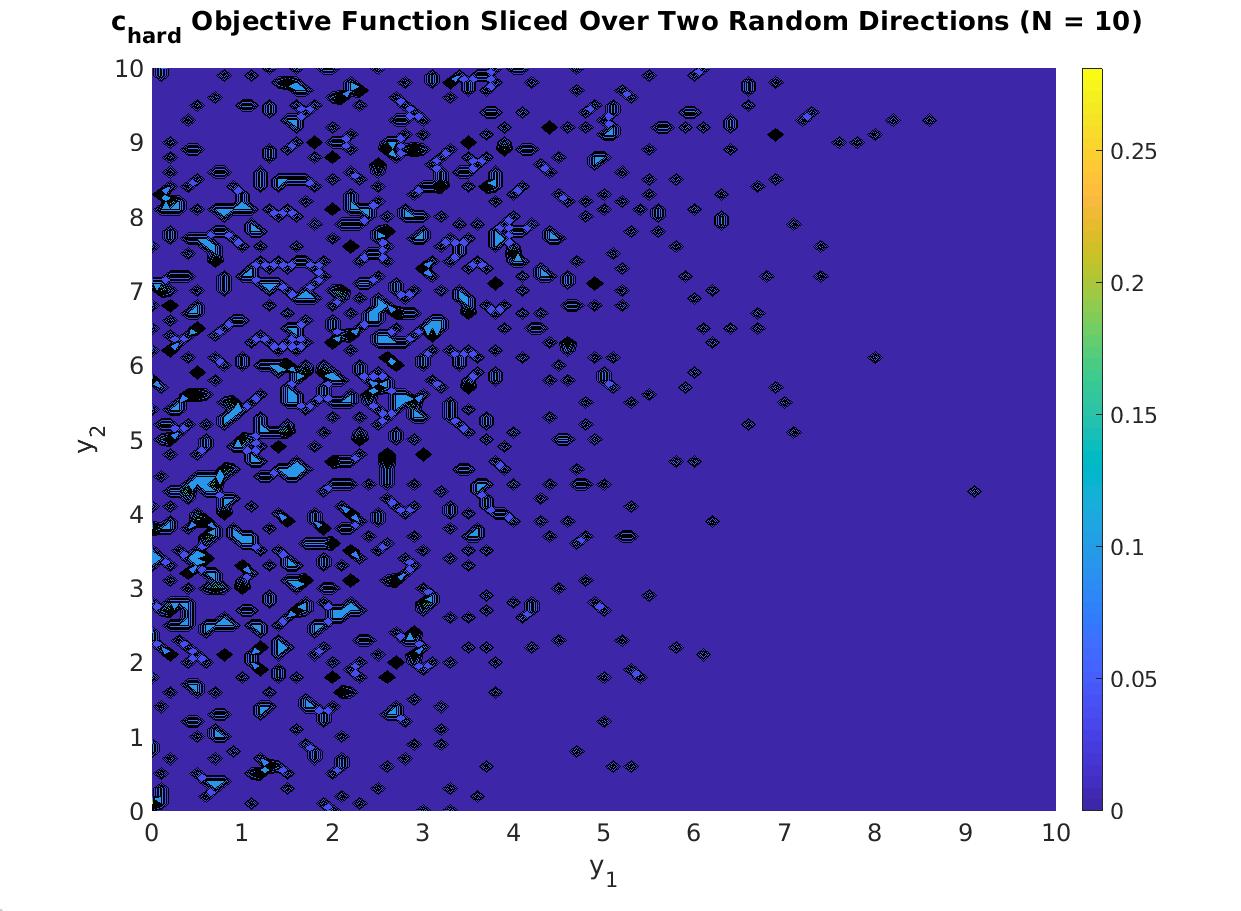}
  & \includegraphics[width=5.0cm]{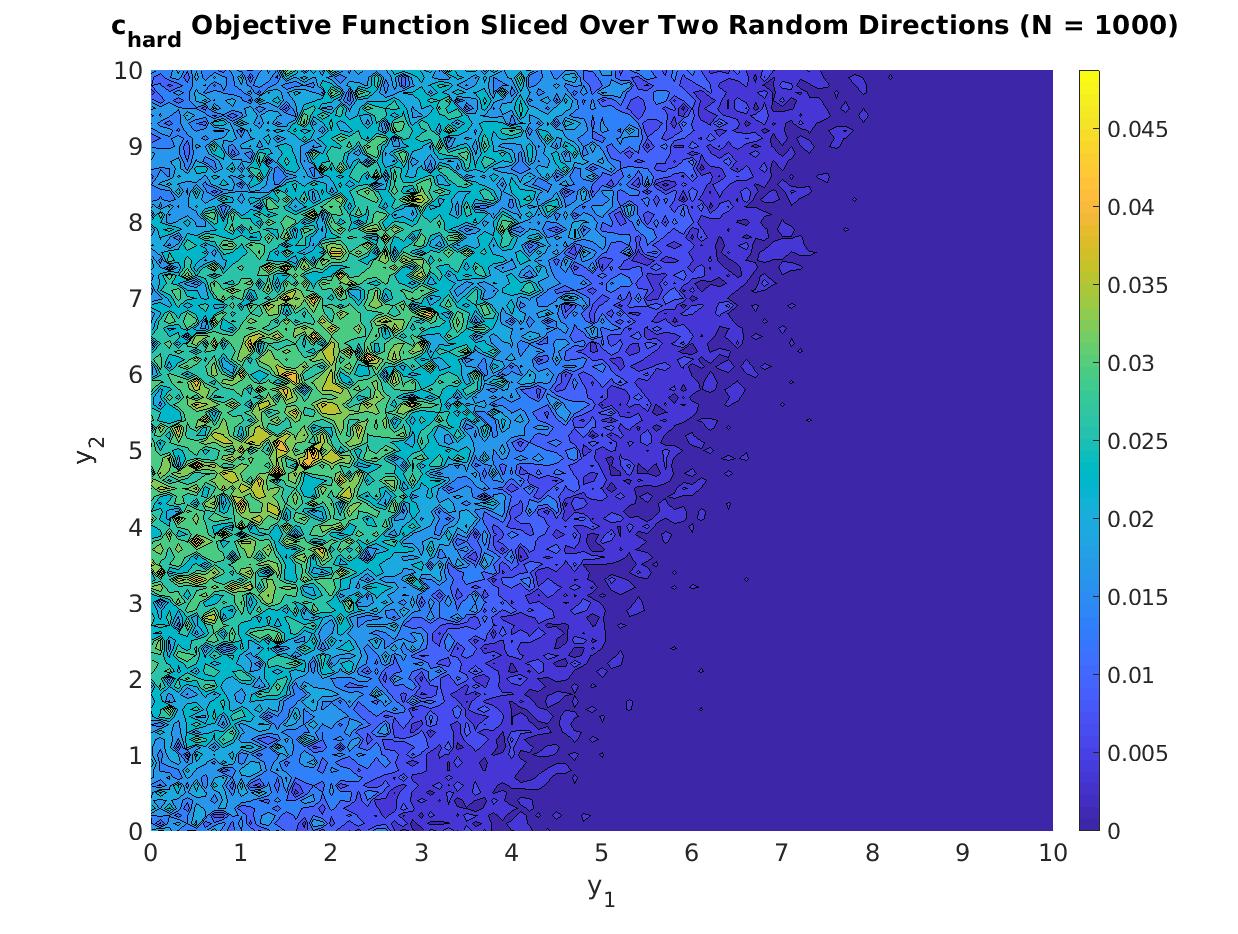}
  \end{tabular}
  \caption{The landscape of the objective function for maximizing the probability of hitting event $c_{hard}$ computed over two random directions $\bfy_1$ and $\bfy_2$. Note the many local maxima, and the objective function's overall non-convexity. Also, note the confusing effects of noise, such as overestimating probabilities when $N$ is small.}
  \label{fig:NorthStar-c_hard-objective-function-landscape}
\end{figure}

In Figure~\ref{fig:NorthStar-c_hard-objective-function-landscape}, we plot the objective function for maximizing $p(c_{hard})$ sliced over two random directions $\bfy_1$ and $\bfy_2$, for $N=10$ and $N=1000$ simulator runs per point respectively. The uniform test template $\bft_{uni}$, defined by
\begin{equation*}
IW = (0.2, 0.2, 0.2, 0.2, 0.2)
\end{equation*}
\begin{equation*}
SW=TW=(0.125, 0.125, 0.125, 0.125, 0.125, 0.125, 0.125, 0.125)
\end{equation*}
\begin{equation*}
CW=(0.5, 0.5)
\end{equation*}
defines the origin in Figure~\ref{fig:NorthStar-c_hard-objective-function-landscape}, and is the starting point for all the optimization experiments in the following subsections. Once again, evaluating the objective function first consists of passing a vector in $\mathbb{R}^{23}$ through standard softmax functions over the $1^{st}$ to $5^{th}$, $6^{th}$ to $13^{th}$, $14^{th}$ to $21^{st}$, and $22^{nd}$ to $23^{rd}$ components. After using the standard softmax function to ensure the 4 directive weight vectors define valid probability distributions, we pass the template defined by these 4 directive weight vectors to the biased random stimuli generator, and track the coverage of the generated random programs over 100 cycles. Note the many local minima, the objective function's overall non-convexity, and how increasing the number of simulator runs per point does not substantially alleviate the non-convexity.

\subsection{Implicit Filtering Technique}
\label{sec:implicit-filtering-numerical-experiments}
\begin{figure*}
  \centering
  \includegraphics[width=5in]{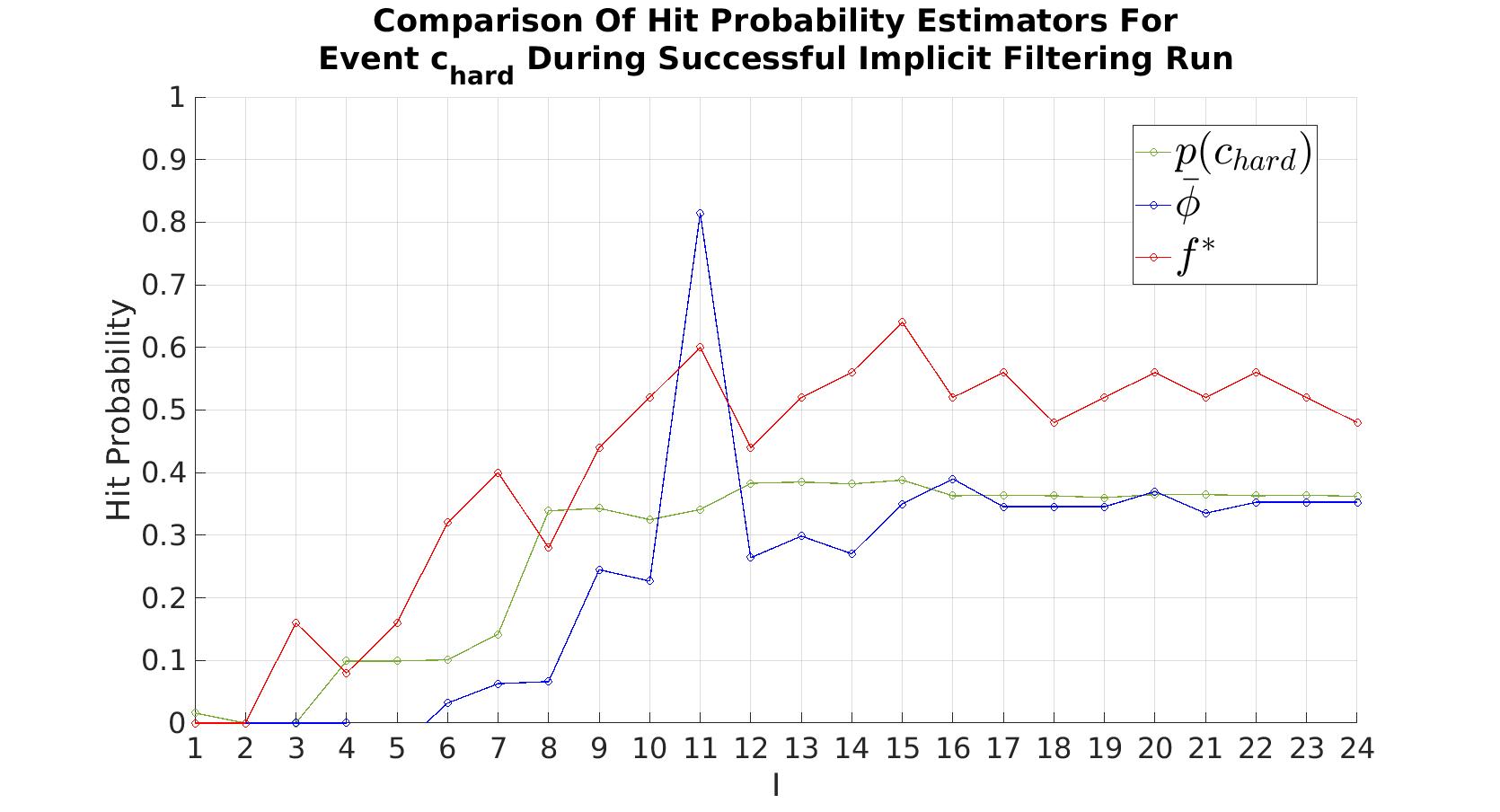}
  \caption{Visualizing the evolution of the hit probability estimators from Table~\ref{tab:implicit_filtering_success} during the successful implicit filtering run. $p(c_{hard})$ denotes the reference ``true" value, which is calculated by averaging over 100000 runs of the simulator at a given template $\bft$. The x-axis $I$ is the iteration number. Note how $f^{*}$ consistently overestimates the probability, whereas $\bar \phi$ generally understimates the probability, but then converges to the ``true" value as the algorithm progresses.}
  \label{fig:NIF_success_hit_probability_estimates}
\end{figure*}

We use implicit filtering to maximize $p(c_{hard})$, which is equivalent to minimizing $-p(c_{hard})$. Table~\ref{tab:implicit_filtering_success} shows the results of a typical successful run of the implicit filtering algorithm. $N$ denotes the number of simulator runs we use to estimate $\bfe_N(\bft)$ at each template, and $n$ is the number of random directions $\bfv$. The algorithm was set to terminate after 50 iterations, or the stencil size $h$ decreased below $1\text{e-}3$. Overall, with a modest budget of 15000 total simulations\footnote{(24 Iterations) $\times (25 \frac{\text{Points}}{\text{Iteration}}) \times (25 \frac{\text{Simulations}}{\text{Point}}) = 15000$ Simulations}, we are able to automatically get within $0.01$ of the best hit probability of $p(c_{hard})=0.4$. Additionally,  we also present the value of $\bar \phi$ estimated by fitting the regularized linear model defined by equation~\ref{eqn:regression-objective} at each iteration. Figure~\ref{fig:NIF_success_hit_probability_estimates} compares the behaviour of $f^{*}$, $\bar \phi$, and the ``true" value of $p(c_{hard})$ averaged over 100000 simulator runs. It is worth noting that the fitted $\bar \phi$ value appears to be a better estimator of $p(c_{hard})$ than the values of $f^{*}$, which is not unexpected given it incorporates information from more runs of the simulator, and nearby points.

\begin{table}[htbp]
\begin{center}
  {\def\arraystretch{1.5}\tabcolsep=6pt
  \begin{tabular}{ | c | c | c | c | c | c | }
    \hline
    \multicolumn{6}{ | c | }{ \textbf{Implicit Filtering History} ($N = 25$, $n=25$, $h_{init} = 50$) } \\
    \hline
    \textit{I} & \textit{$f^{*}$} & \textit{$\bar \phi$} & \textit{Update $\bft_{opt}$?} & \textit{$h$} & $p(c_{hard})$ \\
    \hline
    1 & 0 & 0 & True & 50 & 0.016 \\
    \hline
    2 & 0 & 0 & False & 50 & 0 \\
    \hline
    3 & 0.160 & 7.15\text{e-}5 & True & 25 & 0 \\
    \hline
    4 & 0.080 & 8.62\text{e-}5 & False & 25 & 0.099 \\
    \hline
    5 & 0.160 & 0.043 & False & 12.5 & 0.099 \\
    \hline
    6 & 0.320 & 0.032 & True & 6.25 & 0.101 \\
    \hline
    7 & 0.400 & 0.063 & True & 6.25 & 0.142 \\
    \hline
    8 & 0.280 & 0.066 & False & 6.25 & 0.339 \\
    \hline
    9 & 0.440 & 0.245 & True & 3.125 & 0.343 \\
    \hline
    10 & 0.520 & 0.227 & True & 3.125 & 0.325 \\
    \hline 
    11 & 0.600 & 0.814 & True & 3.125 & 0.341 \\
    \hline
    12 & 0.440 & 0.264 & False & 3.125 & 0.383 \\
    \hline
    13 & 0.520 & 0.297 & False & 1.5625 & 0.385 \\
    \hline
    14 & 0.560 & 0.270 & False & 7.8125\text{e-}1 & 0.382 \\
    \hline
    15 & 0.640 & 0.350 & True & 3.90625\text{e-}1 & 0.388 \\
    \hline
    16 & 0.520 & 0.390 & False & 3.90625\text{e-}1 & 0.363 \\
    \hline
    17 & 0.560 & 0.346 & False & 1.953125\text{e-}1 & 0.364 \\
    \hline
    18 & 0.480 & 0.346 & False & 9.765625\text{e-}2 & 0.363  \\
    \hline
    19 & 0.520 & 0.346 & False & 4.8828125\text{e-}2 & 0.360 \\
    \hline
    20 & 0.560 & 0.370 & False & 2.44140625\text{e-}2 & 0.365 \\
    \hline
    21 & 0.520 & 0.335 & False & 1.220703125\text{e-}2 & 0.365 \\
    \hline
    22 & 0.560 & 0.353 & False & 6.103515625\text{e-}3 & 0.363 \\
    \hline
    23 & 0.520 & 0.353 & False & 3.0517578125\text{e-}3 & 0.364 \\
    \hline
    24 & 0.480 & 0.353 & False & 1.52587890625\text{e-}3 & 0.362 \\
    \hline
    \multicolumn{6}{ | c | }{  \textbf{Summary of Final Results} } \\
    \hline
    \multicolumn{6}{ | c | }{ Total \# of Simulations $= 15000$ } \\
    \hline
    \multicolumn{6}{ | c | }{ $IW_{opt} = [0.5790, 0.2010, 0, 0.2151, 0.0049]$ } \\
    \hline
    \multicolumn{6}{ | c | }{ $SW_{opt} = [0,1,0,0,0,0,0,0]$ } \\
    \hline
    \multicolumn{6}{ | c | }{ $TW_{opt} = [1,0,0,0,0,0,0,0]$ } \\
    \hline
    \multicolumn{6}{ | c | }{ $CW_{opt} = [1,0]$ } \\
    \hline
    \multicolumn{6}{ | c | }{ $f_{opt} = 0.64$, $\bar \phi_{opt} = 0.35$, $p_{opt}(c_{hard}) = 0.39$ } \\
    \hline

  \end{tabular}
  }
\end{center}
\caption{Summary of a successful run of the implicit filtering algorithm. The algorithm is able to \emph{automatically} get within $0.01$ of the best $c_{hard}$ hit probability of $0.4$ using a modest budget of $15000$ total simulations. The optimization was initialized at the uniform template $\bft_{uni}$, and was set to terminate after 50 iterations were exceeded, or the stencil size shrank below $1\text{e-}3$. $I$ is the iteration number.}
\label{tab:implicit_filtering_success}
\end{table}

Table~\ref{tab:implicit_filtering_success} only shows a typical successful run of implicit filtering. The algorithm can also unsuccessfully terminate yielding templates achieving $p_{opt}(c_{hard}) = 0$. As a result, we experiment with the expected performance of the algorithm for different parameter values in Table~\ref{tab:implicit_filtering_ensemble_results}. We also investigate the tradeoff between the number of samples $N$ used to estimate $\bfe_N(\bft)$ at each template, and the number of random directions $n$, at each iteration. For each set of $N$ and $n$ values, we ensemble results over 25 independent runs, and, as in Table~\ref{tab:implicit_filtering_success}, $h_{init}=50$, and the algorithm was set to terminate after 50 iterations, or the stencil size decreased below $1\text{e-}3$.

\begin{table}[htbp]
\begin{center}
  {\def\arraystretch{1.5}\tabcolsep=3pt
  \begin{tabular}{ | c | c | c | c | c | c | c | c | c | c | }
    \hline
    \multicolumn{10}{ | c | }{ \textbf{Per Iteration Budget} = 100 } \\
    \hline
    N & n & $\overline{I}$ & $s^{2}[I]$ & $\overline{f_{opt}}$ & $s^{2}[f_{opt}]$ & $\overline{p_{opt}(c_{hard})}$ & $s^{2}[p_{opt}(c_{hard})]$ & $\max{\{p_{opt}(c_{hard})\}}$ & Failures \\ 
    \hline
    5 & 20 & 17.2 & 0.4 & 0.040 & 0.027 & 0.015 & 0.005 & 0.367 & $\frac{23}{25}$ \\
    \hline
    10 & 10 & 17.2 & 0.8 & 0.048 & 0.028 & 0.022 & 0.006 & 0.335 & $\frac{23}{25}$ \\
    \hline
    20 & 5 & 17.4 & 1.8 & 0.062 & 0.031 & 0.033 & 0.009 & 0.352 & $\frac{22}{25}$ \\
    \hline
    \multicolumn{10}{ | c | }{ \textbf{Per Iteration Budget} = 625 } \\
    \hline
    N & n & $\overline{I}$ & $s^{2}[I]$ & $\overline{f_{opt}}$ & $s^{2}[f_{opt}]$ & $\overline{p_{opt}(c_{hard})}$ & $s^{2}[p_{opt}(c_{hard})]$ & $\max{\{p_{opt}(c_{hard})\}}$ & Failures \\ 
    \hline
    5 & 125 & 18.3 & 3.8 & 0.344 & 0.225 & 0.104 & 0.022 & 0.365 & $\frac{16}{25}$ \\
    \hline
    25 & 25 & 18.5 & 5.7 & 0.192 & 0.083 & 0.104 & 0.025 & 0.377 & $\frac{17}{25}$ \\
    \hline
    125 & 5 & 17.3 & 2.0 & 0.017 & 0.007 & 0.014 & 0.005 & 0.337 & $\frac{24}{25}$ \\
    \hline
    \multicolumn{10}{ | c | }{ \textbf{Per Iteration Budget} = 1250 } \\
    \hline
    N & n & $\overline{I}$ & $s^{2}[I]$ & $\overline{f_{opt}}$ & $s^{2}[f_{opt}]$ & $\overline{p_{opt}(c_{hard})}$ & $s^{2}[p_{opt}(c_{hard})]$ & $\max{\{p_{opt}(c_{hard})\}}$ & Failures \\ 
    \hline
    5 & 250 & 18.1 & 2.7 & 0.320 & 0.227 & 0.105 & 0.025 & 0.375 & $\frac{17}{25}$ \\
    \hline
    10 & 125 & 19.4 & 5.1 & 0.444 & 0.164 & 0.183 & 0.029 & 0.396 & $\frac{11}{25}$ \\
    \hline
    25 & 50 & 19.5 & 10.4 & 0.259 & 0.106 & 0.143 & 0.032 & 0.394 & $\frac{15}{25}$ \\
    \hline
    50 & 25 & 17.6 & 2.8 & 0.066 & 0.033 & 0.041 & 0.013 & 0.388 & $\frac{22}{25}$ \\
    \hline
    125 & 10 & 17.7 & 4.5 & 0.058 & 0.026 & 0.044 & 0.015 & 0.394 & $\frac{22}{25}$ \\
    \hline
    250 & 5 & 17.4 & 3.2 & 0.014 & 0.005 & 0.012 & 0.004 & 0.299 & $\frac{24}{25}$ \\
    \hline
    \multicolumn{10}{ | c | }{ \textbf{Per Iteration Budget} = 2500 } \\
    \hline
    N & n & $\overline{I}$ & $s^{2}[I]$ & $\overline{f_{opt}}$ & $s^{2}[f_{opt}]$ & $\overline{p_{opt}(c_{hard})}$ & $s^{2}[p_{opt}(c_{hard})]$ & $\max{\{p_{opt}(c_{hard})\}}$ & Failures \\ 
    \hline
    5 & 500 & 18.5 & 2.0 & 0.600 & 0.250 & 0.207 & 0.031 & 0.395 & $\frac{10}{25}$ \\
    \hline
    10 & 250 & 20.0 & 5.3 & 0.576 & 0.166 & 0.230 & 0.027 & 0.390 & $\frac{8}{25}$ \\
    \hline
    25 & 100 & 20.5 & 11.4 & 0.379 & 0.119 & 0.201 & 0.034 & 0.390 & $\frac{11}{25}$ \\
    \hline
    50 & 50 & 18.9 & 10.0 & 0.154 & 0.064 & 0.101 & 0.027 & 0.376 & $\frac{18}{25}$ \\
    \hline
    100 & 25 & 18.5 & 9.8 & 0.102 & 0.043 & 0.073 & 0.022 & 0.386 & $\frac{20}{25}$ \\
    \hline
    250 & 10 & 17.6 & 3.9 & 0.036 & 0.016 & 0.031 & 0.011 & 0.390 & $\frac{23}{25}$ \\
    \hline
    500 & 5 & 17.4 & 2.3 & 0.035 & 0.014 & 0.029 & 0.010 & 0.375 & $\frac{23}{25}$ \\
    \hline
    \multicolumn{10}{ | c | }{ \textbf{Per Iteration Budget} = 5000 } \\
    \hline
    N & n & $\overline{I}$ & $s^{2}[I]$ & $\overline{f_{opt}}$ & $s^{2}[f_{opt}]$ & $\overline{p_{opt}(c_{hard})}$ & $s^{2}[p_{opt}(c_{hard})]$ & $\max{\{p_{opt}(c_{hard})\}}$ & Failures \\ 
    \hline
    5 & 1000 & 19.7 & 1.1 & 0.912 & 0.077 & 0.281 & 0.010 & 0.389 & $\frac{2}{25}$ \\
    \hline
    10 & 500 & 20.8 & 4.0 & 0.748 & 0.113 & 0.285 & 0.018 & 0.395 & $\frac{4}{25}$ \\
    \hline
    25 & 200 & 20.8 & 5.8 & 0.539 & 0.080 & 0.276 & 0.021 & 0.393 & $\frac{5}{25}$ \\
    \hline
    50 & 100 & 21.3 & 12.5 & 0.404 & 0.081 & 0.246 & 0.030 & 0.389 & $\frac{8}{25}$ \\
    \hline
    100 & 50 & 20.7 & 13.7 & 0.294 & 0.072 & 0.204 & 0.034 & 0.392 & $\frac{11}{25}$ \\
    \hline
    200 & 25 & 19.0 & 12.8 & 0.131 & 0.046 & 0.103 & 0.029 & 0.394 & $\frac{18}{25}$ \\
    \hline
    500 & 10 & 18.6 & 11.5 & 0.084 & 0.030 & 0.074 & 0.023 & 0.397 & $\frac{20}{25}$ \\
    \hline
    1000 & 5 & 18.1 & 7.0 & 0.048 & 0.017 & 0.043 & 0.014 & 0.392 & $\frac{21}{25}$ \\
    \hline

  \end{tabular}
  }
\end{center}
\caption{Results of the implicit filtering based optimization technique compared over fixed per iteration budgets. Statistics are calculated over $25$ independent runs for each combination of $n$ and $N$, where a bar represents the sample average, and $s^{2}[\cdot]$ represents the sample variance. $N$ is equivalent to simulator runs, and $n$ is the number of directions in which we choose new test templates $\bft$. $I$ is the number of iterations to termination, which occurs when the stencil size parameter reaches less than $h = 0.001$. A failure occurs when the algorithm terminates at a template with a ``true" probability $p_{opt}(c_{hard})=0$. }
\label{tab:implicit_filtering_ensemble_results}
\end{table}

Overall, we see that the implicit filtering technique is not always very reliable. This is exemplified by the ``Failures" column in Table~\ref{tab:implicit_filtering_ensemble_results}, which shows that even for relatively large per iteration budgets, the algorithm can still fail to {\em ever} hit the event $c_{hard}$. Specifically, a failure is defined as the algorithm terminating at a test template that has a $p_{opt}(c_{hard})=0$, or in other words, even after averaging over 100000 simulations at that template, event $c_{hard}$ was {\em never} hit. As expected, Table~\ref{tab:implicit_filtering_ensemble_results} shows the chances of a failure happening are reduced when the per iteration budget is increased, and the number of random directions $n$ is increased. As a general trend, trading off simulation runs $N$ for random directions $n$, given a fixed per iteration budget, is beneficial for the performance of implicit filtering. Only for very small values of $N$, such as $N=5$, does this trend appear to break down. As it is undesirable for how many different parameter choices  the algorithm fails more than half the time. We now investigate if  a gradient-based algorithm performs better overall.

\subsection{Steepest Descent Technique}
\label{sec:steepest-descent-filtering-numerical-experiments}
Now, we use Algorithm~\ref{alg:gradient-based-steepest-descent} to maximize $p(c_{hard})$, which is again equivalent to minimizing $-p(c_{hard})$. Table~\ref{tab:steepest_descent_success} shows the results of a typical successful run of the steepest descent algorithm. The algorithm was set to terminate after 50 iterations, or the line search break flag was set after 10 consecutive line search failures. The line search parameter $\mu_{ls}$ was initialized to $10$. Overall, with a modest budget of 21875 total simulations\footnote{(35 Total Iterations) $\times (25 \frac{\text{Points}}{\text{Iteration}}) \times (25 \frac{\text{Simulations}}{\text{Point}}) = 21875$ Simulations}, we are able to automatically get within $0.09$ of the best hit probability of $p(c_{hard})=0.4$. Note that the term ``total iterations" refers to all the iterations requiring computations, including the failed line searches. For example, iteration 12 in Table~\ref{tab:steepest_descent_success} contributed 4 total iterations, as iteration 12 required 4 line search iterations.  

\begin{table}[htbp]
\begin{center}
  {\def\arraystretch{1.5}\tabcolsep=6pt
  \begin{tabular}{ | c | c | c | c | c | c | c | }
    \hline
    \multicolumn{7}{ | c | }{ \textbf{Steepest Descent History} ($N = 25$, $n=25$, $h = 5$) } \\
    \hline
    \textit{I} & \textit{$\bar \phi$} & \textit{$\| \bfg \|$} & $\| \omega \|$ & \textit{$\mu_{ls}$} & \textit{Update $\bft_{opt}$?} & $p(c_{hard})$ \\
    \hline
    1.1 & 0.018 & 0.0358 & 1.24\text{e-}3 & 10 & True & 0.017 \\
    \hline
    2.1 & 0.029 & 0.1149 & 1.75\text{e-}4 & 20 & True & 0.022 \\
    \hline
    3.4 & 0.030 & 0.0483 & 2.22\text{e-}4 & 5 & True & 0.022 \\
    \hline
    4.1 & 0.027 & 0.0336 & 3.37\text{e-}3 & 5 & False & 0.022 \\
    \hline
    5.1 & 0.027 & 0.0336 & 7.73\text{e-}3 & 10 & False & 0.025 \\
    \hline
    6.1 & 0.021 & 0.0231 & 5.10\text{e-}3 & 20 & False & 0.028 \\
    \hline
    7.1 & 0.021 & 0.0218 & 7.31\text{e-}3 & 40 & False & 0.040 \\
    \hline
    8.1 & 0.026 & 0.0249 & 8.59\text{e-}3 & 80 & False & 0.062 \\
    \hline
    9.1 & 0.030 & 0.0221 & 1.01\text{e-}2 & 160 & False & 0.088 \\
    \hline
    10.1 & 0.034 & 0.0452 & 3.49\text{e-}3 & 320 & True & 0.027 \\
    \hline 
    11.3 & 0.145 & 0.1747 & 4.68\text{e-}3 & 160 & True & 0.333 \\
    \hline
    12.4 & 0.147 & 0.1313 & 3.49\text{e-}2 & 20 & True & 0.339 \\
    \hline
    13.1 & 0.155 & 0.1118 & 1.70\text{e-}2 & 20 & True & 0.331 \\
    \hline
    14.3 & 0.241 & 0.2345 & 2.20\text{e-}4 & 10 & True & 0.312 \\
    \hline
    15.10 & 0.031 & 0.2225 & 5.89\text{e-}4 & 1.953125\text{e-}2 & False & 0.314 \\
    \hline
    \multicolumn{7}{ | c | }{  \textbf{Summary of Final Results} } \\
    \hline
    \multicolumn{7}{ | c | }{ Total \# of Simulations $= 21875$ } \\
    \hline
    \multicolumn{7}{ | c | }{ $IW_{opt} = [0.4377, 0.1642, 0.2597, 0.1368, 0.0015]$ } \\
    \hline
    \multicolumn{7}{ | c | }{ $SW_{opt} = [0,0,0,0,1,0,0,0]$ } \\
    \hline
    \multicolumn{7}{ | c | }{ $TW_{opt} = [0,1,0,0,0,0,0,0]$ } \\
    \hline
    \multicolumn{7}{ | c | }{ $CW_{opt} = [1,0]$ } \\
    \hline
    \multicolumn{7}{ | c | }{ $\bar \phi_{opt} = 0.24$, $p_{opt}(c_{hard}) = 0.31$ } \\
    \hline

  \end{tabular}
  }
\end{center}
\caption{Summary of a successful run of the steepest descent algorithm (Algorithm~\ref{alg:gradient-based-steepest-descent}). The algorithm is able to \emph{automatically} get within $0.09$ of the best $c_{hard}$ hit probability of $0.4$ using a modest budget of $21875$ total simulations. The optimization was initialized at the uniform template $\bft_{uni}$, and was set to terminate after 50 main iterations were exceeded, or after 10 consecutive line search failures. $I$ is the iteration number, formatted so the number after the decimal point represents the final line search iteration for the given main iteration.}
\label{tab:steepest_descent_success}
\end{table}

\begin{figure*}
  \centering
  \includegraphics[width=5in]{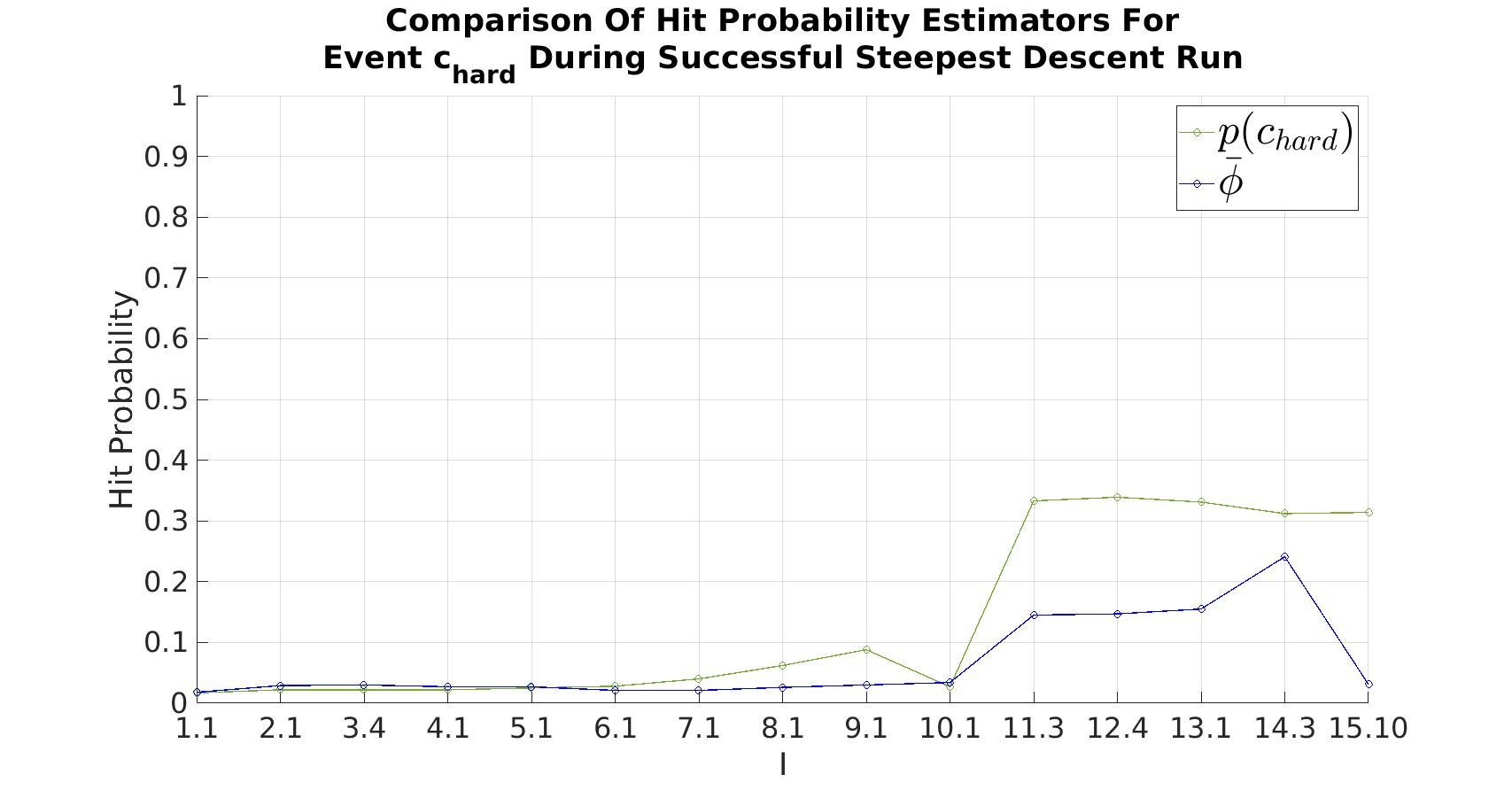}
  \caption{Visualizing the evolution of the hit probability estimators from Table~\ref{tab:steepest_descent_success} during the successful run of the steepest descent algorithm (Algorithm~\ref{alg:gradient-based-steepest-descent}). $p(c_{hard})$ denotes the reference ``true" value, which is calculated by averaging over 100000 runs of the simulator at a given template $\bft$. The x-axis shows the iteration number, where the value after the decimal point is the final line search iteration number for the given main iteration. Note how $\bar \phi$ understimates the ``true" probability, especially towards the end of the run.}
  \label{fig:SD_success_hit_probability_estimates}
\end{figure*}

Figure~\ref{fig:SD_success_hit_probability_estimates} compares the behaviour of $\bar \phi$, and the ``true" value of $p(c_{hard})$. In general, $\bar \phi$ tracks $p(c_{hard})$ closely, but with a tendency to vary more slowly. This is because of the averaging effect of the algorithm. %Given the strong non-convexity of the objective function, exemplified by Figure~\ref{fig:NorthStar-c_hard-objective-function-landscape}, it is likely beneficial to adapt the stencil size $h$ as Algorithm~\ref{alg:gradient-based-steepest-descent} progresses. Further potential evidence for adapting $h$ during Algorithm~\ref{alg:gradient-based-steepest-descent} being beneficial is presented in Table~\ref{tab:steepest_descent_ensemble_results}. However, we do not explore adapting the size of $h$ during Algorithm~\ref{alg:gradient-based-steepest-descent} further in this paper.

Similar to Table~\ref{tab:implicit_filtering_ensemble_results}, Table~\ref{tab:steepest_descent_ensemble_results} investigates the expected performance of Algorithm~\ref{alg:gradient-based-steepest-descent} for different parameter values. Like with Table~\ref{tab:implicit_filtering_ensemble_results}, for each set of $N$ and $n$ values, we ensemble over 25 independent runs, and, as in Table~\ref{tab:steepest_descent_success}, $h=5$ and $\mu_{ls}$ is initialized to $10$, and the algorithm was set to terminate after 50 iterations, or the line search break flag was set after 10 consecutive line search failures.

\begin{table}[htbp]
\begin{center}
  {\def\arraystretch{1.5}\tabcolsep=3pt
  \begin{tabular}{ | c | c | c | c | c | c | c | c | c | c | }
    \hline
    \multicolumn{10}{ | c | }{ \textbf{Per Iteration Budget} = 100 } \\
    \hline
    N & n & $\overline{I_{t}}$ & $s^{2}[I_{t}]$ & $\overline{\bar \phi_{opt}}$ & $s^{2}[\bar \phi_{opt}]$ & $\overline{p_{opt}(c_{hard})}$ & $s^{2}[p_{opt}(c_{hard})]$ & $\max{\{p_{opt}(c_{hard})\}}$ & Failures \\
    \hline
    5 & 20 & 82.5 & 290.2 & 0.164 & 0.013 & 0.124 & 0.016 & 0.375 & $\frac{0}{25}$ \\
    \hline
    10 & 10 & 81.7 & 1.04\text{e}3 & 0.041 & 7.13\text{e-}4 & 0.077 & 0.010 & 0.372 &  $\frac{3}{25}$ \\
    \hline
    20 & 5 & 93.8 & 168.6 & 0.020 & 2.41\text{e-}4 & 0.022 & 0.002 & 0.165 & $\frac{12}{25}$ \\
    \hline
    \multicolumn{10}{ | c | }{ \textbf{Per Iteration Budget} = 625 } \\
    \hline
    N & n & $\overline{I_{t}}$ & $s^{2}[I_{t}]$ & $\overline{\bar \phi_{opt}}$ & $s^{2}[\bar \phi_{opt}]$ & $\overline{p_{opt}(c_{hard})}$ & $s^{2}[p_{opt}(c_{hard})]$ & $\max{\{p_{opt}(c_{hard})\}}$ & Failures \\
    \hline
    5 & 125 & 79.6 & 2.1 & 0.180 & 7.27\text{e-}5 & 0.355 & 0.001 & 0.397 & $\frac{0}{25}$ \\
    \hline
    25 & 25 & 47.7 & 645.9 & 0.337 & 0.085 & 0.163 & 0.017 & 0.367 & $\frac{0}{25}$ \\
    \hline
    125 & 5 & 94 & 246.0 & 0.020 & 1.47\text{e-}4 & 0.018 & 0.001 & 0.159 & $\frac{11}{25}$ \\
    \hline
    \multicolumn{10}{ | c | }{ \textbf{Per Iteration Budget} = 1250 } \\
    \hline
    N & n & $\overline{I_{t}}$ & $s^{2}[I_{t}]$ & $\overline{\bar \phi_{opt}}$ & $s^{2}[\bar \phi_{opt}]$ & $\overline{p_{opt}(c_{hard})}$ & $s^{2}[p_{opt}(c_{hard})]$ & $\max{\{p_{opt}(c_{hard})\}}$ & Failures \\
    \hline
    5 & 250 & 79.7 & 2.0 & 0.183 & 3.48\text{e-}5 & 0.375 & 3.44\text{e-}4 & 0.399 & $\frac{0}{25}$ \\
    \hline
    10 & 125 & 80.3 & 2.3 & 0.182 & 4.00\text{e-}5 & 0.362 & 8.63\text{e-}4 & 0.400 & $\frac{0}{25}$ \\
    \hline
    25 & 50 & 80.5 & 3.0 & 0.181 & 2.07\text{e-}4 & 0.340 & 0.001 & 0.384 & $\frac{0}{25}$ \\
    \hline
    50 & 25 & 43.5 & 646.3 & 0.211 & 0.024 & 0.151 & 0.019 & 0.394 & $\frac{0}{25}$ \\
    \hline
    125 & 10 & 71.2 & 951.2 & 0.061 & 0.003 & 0.118 & 0.011 & 0.379 & $\frac{0}{25}$ \\
    \hline
    250 & 5 & 98.0 & 336.5 & 0.019 & 8.95\text{e-}5 & 0.036 & 0.004 & 0.218 & $\frac{10}{25}$\\
    \hline
    \multicolumn{10}{ | c | }{ \textbf{Per Iteration Budget} = 2500 } \\
    \hline
    N & n & $\overline{I_{t}}$ & $s^{2}[I_{t}]$ & $\overline{\bar \phi_{opt}}$ & $s^{2}[\bar \phi_{opt}]$ & $\overline{p_{opt}(c_{hard})}$ & $s^{2}[p_{opt}(c_{hard})]$ & $\max{\{p_{opt}(c_{hard})\}}$ & Failures \\
    \hline
    5 & 500 & 80.2 & 1.8 & 0.185 & 1.40\text{e-}5 & 0.381 & 2.15\text{e-}4 & 0.402 & $\frac{0}{25}$ \\
    \hline
    10 & 250 & 79.6 & 1.8 & 0.185 & 2.18\text{e-}5 & 0.381 & 2.00\text{e-}4 & 0.400 & $\frac{0}{25}$ \\
    \hline
    25 & 100 & 79.2 & 2.5 & 0.184 & 3.33\text{e-}5 & 0.358 & 0.001 & 0.397 & $\frac{0}{25}$ \\
    \hline
    50 & 50 & 79.4 & 1.8 & 0.183 & 1.04\text{e-}4 & 0.338 & 0.002 & 0.398 & $\frac{0}{25}$ \\
    \hline
    100 & 25 & 44.2 & 585.1 & 0.250 & 0.043 & 0.189 & 0.017 & 0.382 & $\frac{0}{25}$ \\
    \hline
    250 & 10 & 80.0 & 1.76\text{e}3 & 0.054 & 0.002 & 0.120 & 0.012 & 0.344 & $\frac{0}{25}$ \\
    \hline
    500 & 5 & 101.5 & 99.0 & 0.019 & 8.10\text{e-}5 & 0.044 & 0.006 & 0.331 & $\frac{7}{25}$ \\
    \hline
    \multicolumn{10}{ | c | }{ \textbf{Per Iteration Budget} = 5000 } \\
    \hline
    N & n & $\overline{I_{t}}$ & $s^{2}[I_{t}]$ & $\overline{\bar \phi_{opt}}$ & $s^{2}[\bar \phi_{opt}]$ & $\overline{p_{opt}(c_{hard})}$ & $s^{2}[p_{opt}(c_{hard})]$ & $\max{\{p_{opt}(c_{hard})\}}$ & Failures \\
    \hline
    5 & 1000 & 80.6 & 1.6 & 0.185 & 1.01\text{e-}5 & 0.385 & 2.02\text{e-}4 & 0.400 & $\frac{0}{25}$ \\
    \hline
    10 & 500 & 80.1 & 3.3 & 0.185 & 1.16\text{e-}5 & 0.387 & 8.26\text{e-}5 & 0.402 & $\frac{0}{25}$ \\
    \hline
    25 & 200 & 80.0 & 2.4 & 0.186 & 2.28\text{e-}5 & 0.378 & 5.03\text{e-}4 & 0.402 & $\frac{0}{25}$ \\
    \hline
    50 & 100 & 80.0 & 2.1 & 0.182 & 4.65\text{e-}5 & 0.364 & 7.69\text{e-}4 & 0.399 & $\frac{0}{25}$ \\
    \hline
    100 & 50 & 79.9 & 1.6 & 0.185 & 1.43\text{e-}4 & 0.354 & 0.002 & 0.395 & $\frac{0}{25}$ \\
    \hline
    200 & 25 & 46.7 & 442.7 & 0.207 & 0.012 & 0.177 & 0.018 & 0.376 & $\frac{0}{25}$ \\
    \hline
    500 & 10 & 89.1 & 1.65\text{e}3 & 0.057 & 0.001 & 0.126 & 0.015 & 0.384 & $\frac{2}{25}$ \\
    \hline
    1000 & 5 & 99.9 & 326.2 & 0.027 & 4.83\text{e-}4 & 0.052 & 0.005 & 0.213 & $\frac{7}{25}$ \\
    \hline

  \end{tabular}
  }
\end{center}
\caption{Results of the steepest descent based optimization technique compared over fixed per iteration budgets. Statistics are calculated over $25$ independent runs for each combination of $n$ and $N$, where a bar represents the sample average, and $s^{2}[\cdot]$ represents the sample variance. $N$ is equivalent to simulator runs, and $n$ is the number of directions in which we choose new test templates $\bft$. $I_{t}$ is the number of total iterations to termination, which includes all line searches. Termination occurs after 50 iterations, or 10 consecutive failed line searches. A failure occurs when the algorithm terminates at a template with a ``true" probability $p(c_{hard})=0$. }
\label{tab:steepest_descent_ensemble_results}
\end{table}

Overall, the gradient based steepest descent technique appears much more reliable than the implicit filtering technique. Algorithm~\ref{alg:gradient-based-steepest-descent} almost always terminates at a template that at least hits the event $c_{hard}$ a minimum of once in 100000 simulation runs. It is also worth noting that in most cases $\bar \phi$ underestimates $p(c_{hard})$, sometimes by a large margin of up to almost 0.2. The authors conjecture this may be due to a relatively large choice of $h$, which is not refined during Algorithm~\ref{alg:gradient-based-steepest-descent}. The effects of a relatively large $h$ should be especially pronounced if the optima are rather sharp, which given domain knowledge, is not unlikely for this problem. However, as with the implicit filtering technique, the steepest descent algorithm's performance strongly benefits from trading off $N$ for $n$, given a fixed per iteration budget. For both algorithms, it appears that in general, coarsely sampling many points is preferable to sampling a few points with high accuracy at each point. 

\subsection{BFGS Technique}
\label{sec:BFGS-numerical-experiments}

Finally, following the framework of Algorithm~\ref{alg:gradient-based-steepest-descent}, we use BFGS to minimize $-p(c_{hard})$. Specifically, we use a limited-memory implementation of the BFGS method, also referred to as L-BFGS. To compute the BFGS directions, we use the L-BFGS two-loop recursion detailed on page 225 of \cite{nw}. We set the initial inverse Hessian approximation to be a scaled version of the identity matrix, where the scaling factor is given by equation 9.6 on page 226 of \cite{nw}. As a result, a back tracking line search starting with $\mu_{ls}=1$ at each iteration, and refining by a factor of two for each line search failure, was employed. However, we set the memory value, $m$, for our L-BFGS implementation to $m=100$, which was almost always greater than the number of iterations before termination. As a result, almost all of the time our L-BFGS implementation was equivalent to a standard BFGS implementation.

Of the three algorithms we tested, the L-BFGS implementation had the most difficulty obtaining test templates achieving close to $p(c_{hard})=0.4$, and is not competitive with either implicit filtering or steepest descent. As with Tables~\ref{tab:implicit_filtering_ensemble_results} and~\ref{tab:steepest_descent_ensemble_results}, Table~\ref{tab:LBFGS_ensemble_results} investigates the expected performance of BFGS for different parameter values. Like with Table~\ref{tab:steepest_descent_ensemble_results}, for each set of $N$ and $n$ values, we ensemble over 25 independent runs, $h=5$, and the algorithm was set to terminate after 50 iterations, or the line search break flag was set after 10 consecutive line search failures.

\begin{table}[htbp]
\begin{center}
  {\def\arraystretch{1.5}\tabcolsep=3pt
  \begin{tabular}{ | c | c | c | c | c | c | c | c | c | c | }
    \hline
    \multicolumn{10}{ | c | }{ \textbf{Per Iteration Budget} = 100 } \\
    \hline
    N & n & $\overline{I_{t}}$ & $s^{2}[I_{t}]$ & $\overline{\bar \phi_{opt}}$ & $s^{2}[\bar \phi_{opt}]$ & $\overline{p_{opt}(c_{hard})}$ & $s^{2}[p_{opt}(c_{hard})]$ & $\max{\{p_{opt}(c_{hard})\}}$ & Failures \\
    \hline
    5 & 20 & 32.7 & 423.9 & 0.032 & 4.27\text{e-}4 & 0.017 & 4.37\text{e-}7 & 0.019 & $\frac{0}{25}$ \\
    \hline
    10 & 10 & 45.8 & 497.4 & 0.023 & 1.51\text{e-}4 & 0.017 & 1.08\text{e-}4 & 0.061 & $\frac{1}{25}$ \\
    \hline
    20 & 5 & 33.2 & 682.2 & 0.012 & 8.69\text{e-}5 & 0.017 & 2.56\text{e-}6 & 0.024 & $\frac{0}{25}$ \\
    \hline
    \multicolumn{10}{ | c | }{ \textbf{Per Iteration Budget} = 625 } \\
    \hline
    N & n & $\overline{I_{t}}$ & $s^{2}[I_{t}]$ & $\overline{\bar \phi_{opt}}$ & $s^{2}[\bar \phi_{opt}]$ & $\overline{p_{opt}(c_{hard})}$ & $s^{2}[p_{opt}(c_{hard})]$ & $\max{\{p_{opt}(c_{hard})\}}$ & Failures \\
    \hline
    5 & 125 & 45.0 & 407.1 & 0.018 & 1.39\text{e-}4 & 0.022 & 4.13\text{e-}4 & 0.119 & $\frac{0}{25}$ \\
    \hline
    25 & 25 & 37.0 & 255.1 & 0.072 & 0.001 & 0.018 & 3.02\text{e-}5 & 0.044 & $\frac{0}{25}$ \\
    \hline
    125 & 5 & 57.5 & 187.3 & 0.016 & 5.15\text{e-}5 & 0.016 & 1.14\text{e-}4 & 0.059 & $\frac{2}{25}$ \\
    \hline
    \multicolumn{10}{ | c | }{ \textbf{Per Iteration Budget} = 1250 } \\
    \hline
    N & n & $\overline{I_{t}}$ & $s^{2}[I_{t}]$ & $\overline{\bar \phi_{opt}}$ & $s^{2}[\bar \phi_{opt}]$ & $\overline{p_{opt}(c_{hard})}$ & $s^{2}[p_{opt}(c_{hard})]$ & $\max{\{p_{opt}(c_{hard})\}}$ & Failures \\
    \hline
    5 & 250 & 51.5 & 232.4 & 0.019 & 1.46\text{e-}4 & 0.025 & 3.79\text{e-}4 & 0.112 & $\frac{0}{25}$ \\
    \hline
    10 & 125 & 51.0 & 244.1 & 0.030 & 0.002 & 0.048 & 0.005 & 0.270 & $\frac{0}{25}$ \\
    \hline
    25 & 50 & 49.0 & 215.7 & 0.028 & 9.34\text{e-}4 & 0.032 & 0.001 & 0.169 & $\frac{0}{25}$ \\
    \hline
    50 & 25 & 35.0 & 361.3 & 0.065 & 0.003 & 0.030 & 0.004 & 0.337 & $\frac{0}{25}$ \\
    \hline
    125 & 10 & 51.7 & 315.7 & 0.018 & 1.54\text{e-}4 & 0.018 & 1.36\text{e-}4 & 0.059 & $\frac{1}{25}$ \\
    \hline
    250 & 5 & 53.9 & 275.1 & 0.018 & 1.75\text{e-}4 & 0.016 & 1.34\text{e-}4 & 0.050 & $\frac{3}{25}$\\
    \hline
    \multicolumn{10}{ | c | }{ \textbf{Per Iteration Budget} = 2500 } \\
    \hline
    N & n & $\overline{I_{t}}$ & $s^{2}[I_{t}]$ & $\overline{\bar \phi_{opt}}$ & $s^{2}[\bar \phi_{opt}]$ & $\overline{p_{opt}(c_{hard})}$ & $s^{2}[p_{opt}(c_{hard})]$ & $\max{\{p_{opt}(c_{hard})\}}$ & Failures \\
    \hline
    5 & 500 & 44.9 & 288.2 & 0.027 & 4.52\text{e-}4 & 0.041 & 0.002 & 0.203 & $\frac{0}{25}$ \\
    \hline
    10 & 250 & 46.5 & 282.4 & 0.026 & 3.87\text{e-}4 & 0.035 & 0.001 & 0.142 & $\frac{0}{25}$ \\
    \hline
    25 & 100 & 52 & 367.7 & 0.019 & 2.37\text{e-}4 & 0.026 & 7.01\text{e-}4 & 0.150 & $\frac{0}{25}$ \\
    \hline
    50 & 50 & 47.6 & 353.2 & 0.025 & 4.45\text{e-}4 & 0.031 & 8.01\text{e-}4 & 0.129 & $\frac{0}{25}$ \\
    \hline
    100 & 25 & 34.4 & 223.3 & 0.045 & 0.001 & 0.026 & 6.83\text{e-}4 & 0.144 & $\frac{0}{25}$ \\
    \hline
    250 & 10 & 51.2 & 423.3 & 0.018 & 3.91\text{e-}5 & 0.017 & 5.79\text{e-}5 & 0.045 & $\frac{1}{25}$ \\
    \hline
    500 & 5 & 58.5 & 135.6 & 0.014 & 5.50\text{e-}5 & 0.015 & 4.59\text{e-}5 & 0.026 & $\frac{2}{25}$ \\
    \hline
    \multicolumn{10}{ | c | }{ \textbf{Per Iteration Budget} = 5000 } \\
    \hline
    N & n & $\overline{I_{t}}$ & $s^{2}[I_{t}]$ & $\overline{\bar \phi_{opt}}$ & $s^{2}[\bar \phi_{opt}]$ & $\overline{p_{opt}(c_{hard})}$ & $s^{2}[p_{opt}(c_{hard})]$ & $\max{\{p_{opt}(c_{hard})\}}$ & Failures \\
    \hline
    5 & 1000 & 47.6 & 531.8 & 0.016 & 3.44\text{e-}6 & 0.022 & 1.98\text{e-}5 & 0.032 & $\frac{0}{25}$ \\
    \hline
    10 & 500 & 46.3 & 313.0 & 0.028 & 0.001 & 0.046 & 0.007 & 0.368 & $\frac{0}{25}$ \\
    \hline
    25 & 200 & 47.1 & 424.9 & 0.023 & 4.27\text{e-}4 & 0.033 & 0.001 & 0.161 & $\frac{0}{25}$ \\
    \hline
    50 & 100 & 47.4 & 391.3 & 0.016 & 2.28\text{e-}5 & 0.022 & 9.76\text{e-}5 & 0.065 & $\frac{0}{25}$ \\
    \hline
    100 & 50 & 47.8 & 452.7 & 0.017 & 8.17\text{e-}5 & 0.025 & 2.77\text{e-}4 & 0.094 & $\frac{0}{25}$ \\
    \hline
    200 & 25 & 31.2 & 215.4 & 0.062 & 0.013 & 0.018 & 1.35\text{e-}6 & 0.021 & $\frac{0}{25}$ \\
    \hline
    500 & 10 & 48.0 & 212.9 & 0.020 & 2.84\text{e-}4 & 0.021 & 1.86\text{e-}4 & 0.071 & $\frac{2}{25}$ \\
    \hline
    1000 & 5 & 58.6 & 345.5 & 0.015 & 2.56\text{e-}5 & 0.017 & 4.55\text{e-}5 & 0.037 & $\frac{1}{25}$ \\
    \hline

  \end{tabular}
  }
\end{center}
\caption{Results of the L-BFGS based optimization technique compared over fixed per iteration budgets. Statistics are calculated over $25$ independent runs for each combination of $n$ and $N$, where a bar represents the sample average, and $s^{2}[\cdot]$ represents the sample variance. $N$ is equivalent to simulator runs, and $n$ is the number of directions in which we choose new test templates $\bft$. $I_{t}$ is the number of total iterations to termination, which includes all line searches. Termination occurs after 50 iterations, or 10 consecutive failed line searches. A failure occurs when the algorithm terminates at a template with a ``true" probability $p(c_{hard})=0$. }
\label{tab:LBFGS_ensemble_results}
\end{table}

Whereas Tables~\ref{tab:implicit_filtering_ensemble_results} and~\ref{tab:steepest_descent_ensemble_results} show that with a budget of 625 simulations per iteration, $n=25$, and $N=25$, the implicit filtering and steepest descent techniques on average achieve $p_{opt}(c_{hard})=0.104$ and $p_{opt}(c_{hard})=0.163$ respectively, Table~\ref{tab:LBFGS_ensemble_results} shows the L-BFGS method only achieves $p_{opt}(c_{hard})=0.018$ on average. However, as with the previous two algorithms, there is still a noticeable benefit from trading off $N$ for $n$, given a fixed per iteration budget, and increasing the per iteration budget can improve performance. The L-BFGS technique also fails much less frequently than the implicit filtering technique. Overall though, the ensemble results suggest this method is inferior to the steepest descent based approach, and that the standard BFGS technique may need further modifications to handle the noise in this problem.

%==========================
\section{Summary and Conclusions}
\label{sec:conclusions-sec}
%==========================

In this paper, we have proposed three algorithms for solving the coverage directed generation problem, all based on the key observation that the problem can be posed as derivative free optimization of a noisy objective function. By applying techniques from statistical parameter estimation and inverse problems, namely the generalized cross validation technique, we are able to generate quality estimates of the gradient of the noisy objective function. With these gradient estimates, we are able to build algorithms  that adapt the steepest descent and BFGS techniques from non-noisy continuous optimization. 

The algorithm based on gradient descent, on average, empirically outperforms a simple, but sometimes surprisingly effective, implicit filtering based approach. Numerical experiments with a high-level software model of part of IBM's NorthStar processor show that both the implicit filtering and steepest descent techniques are economical in terms of the total number of simulations required for them to be effective, and how to best choose parameters given a fixed per iteration budget of simulations. Furthermore, all our algorithms are relatively easily parallelized in practice, as the repeated simulations at a single point $N$ can be carried out in parallel, and this can further be done in parallel for the $n$ points along the random directions, with the only major bottleneck being the work required during the decision to update the next template. 

We suspect that the use of Inverse Problem based techniques for gradient estimation
can be further extended to the evaluation of Hessians and in other contexts where the 
function and gradients are noisy, and this will be investigated in the future.

%\subsection{Subsection title}
%\label{sec:2}
%as required. Don't forget to give each section
%and subsection a unique label (see Sect.~\ref{sec:1}).
%\paragraph{Paragraph headings} Use paragraph headings as needed.
%\begin{equation}
%\end{equation}

%\begin{figure}
%  \includegraphics{example.eps}
%\end{figure}

%\begin{figure*}
%  \includegraphics[width=0.75\textwidth]{example.eps}
%\end{figure*}

%\begin{table}
%\begin{tabular}{lll}
%\hline\noalign{\smallskip}
%\end{tabular}
%\end{table}

\begin{acknowledgements}
EH and BI's work is supported by the Natural Sciences and Engineering Research Council of Canada (NSERC). RG, BS, and AZ's work is supported by IBM.
\end{acknowledgements}

\section*{Conflict of interest}

The authors declare that they have no conflict of interest.

\bibliographystyle{spmpsci}      % mathematics and physical sciences
\bibliography{biblio1}   % name your BibTeX data base

%\begin{thebibliography}{}
% etc
%\end{thebibliography}

\end{document}